\documentclass[12pt]{amsart}
\usepackage{amssymb, enumerate, mdwlist, stmaryrd}

\evensidemargin\paperwidth
\advance\evensidemargin-\textwidth
\advance\oddsidemargin-1in
\evensidemargin\oddsidemargin

\topmargin\paperheight
\advance\topmargin-\textheight
\advance\topmargin-1in

\theoremstyle{plain}
\newtheorem{theorem}{Theorem}[section]
\theoremstyle{plain}

\theoremstyle{plain}
\newtheorem{lemma}[theorem]{Lemma}
\theoremstyle{plain}
\newtheorem{corollary}[theorem]{Corollary}
\theoremstyle{plain}

\newtheorem{question}[theorem]{Question}
\theoremstyle{plain}
\newtheorem{mthm}{Theorem}

\theoremstyle{definition}
\newtheorem{definition}[theorem]{Definition}
\theoremstyle{remark}
\newtheorem{remark}[theorem]{Remark}
\theoremstyle{remark}

\theoremstyle{remark}

\title[Multi-way expanders]
{Multi-way expanders and imprimitive group actions on graphs.}
\author{Masato Mimura} \thanks{Supported in part by the Grant-in-Aid for Young Scientists (B), no.25800033 from the JSPS}
\address{Masato Mimura\\
Mathematical Institute, Tohoku University}
\email{mimura-mas@m.tohoku.ac.jp}
\date{\today}

\begin{document}

\begin{abstract}
For $n\geq 2$, the concept of $n$-way expanders was defined by various researchers. Bigger $n$ gives a weaker notion in general, and $2$-way expanders coincide with expanders in usual sense. Koji Fujiwara asked whether these concepts are equivalent to that of ordinary expanders for all $n$ for a sequence of Cayley graphs. In this paper, we answer his question in the affirmative. Furthermore, we obtain universal inequalities on multi-way isoperimetric constants on any finite connected vertex-transitive graph, and show that gaps between these constants imply the imprimitivity of the group action on the graph. 
\end{abstract}

\keywords{Expanders; multi-way isoperimetric constants; imprimitive actions}

\maketitle

\section{Introduction}
\label{sec:intro}
In this paper, let $n$ represent a natural number at least $2$. We assume that all graphs $\Gamma=(V,E)$ are finite, undirected, regular, and without multiple edges or self-loops. Denote by $d_{\Gamma}$  the regularity of $\Gamma$. For a \textit{Cayley graph} $\Gamma=\mathrm{Cay}(G,S)$, we use the right-multiplication to connect edges in order to have the left-action by graph isomorphisms. We allow the case where $\Gamma$ is disconnected (for Cayley graphs, this amounts to saying that $S$ does not generate the whole $G$). For disjoint subsets $A, B$ of the vertex set $V$, $\partial (A,B)$ denotes the \textit{edge boundary}, that means, the set $\{e=(u,v)\in E:u\in A, v\in B\})$, and $\partial A$ denotes $\partial (A,V\setminus A)$. In addition, $\delta (A,B)$ denotes the \textit{symmetric vertex boundary}, that is, the set $\{u \in A : \exists e=(u,v) \in \partial (A,B)\} \sqcup \{v \in B : \exists  e=(u,v) \in \partial (A,B)\})$, and $\delta A$ denotes $\delta (A,V\setminus A)$. For $l\in \mathbb{N}$, by $\mathfrak{S}_l$, we denote the symmetric group of degree $l$. Let $[l]$ be the set $\{1,2,\ldots ,l\}$.

For $(|V|\geq )n\geq 2$, the following three quantities are defined.
\begin{definition}\label{defn:hn}
Let $\Gamma=(V,E)$, and $2\leq n\leq |V|$. 
\begin{enumerate}[$(1)$]
  \item  The $n$\textit{-way isoperimetric constant} $h_n(\Gamma)$ is defined by\[
h_n(\Gamma)=\min \max_{1\leq i\leq n}\frac{|\partial A_i|}{|A_i|}.
\] 
Here the minimum is taken over all partitions of $V$ into $n$ non-empty disjoint subsets $V=\bigsqcup_{i=1}^n A_i$. 
  \item The $n$\textit{-way symmetric vertex isoperimetric constant} $\iota_n(\Gamma)$ is defined by 
\[
\iota_n(\Gamma)=\min \max_{1\leq i\leq n}\frac{|\delta A_i|}{|A_i|}.
\]
Here $(A_1,\ldots ,A_n)$ runs over the same partitions as in item $(1)$.
  \item The $\lambda_n(\Gamma)$ is the $n$\textit{-th nonnegative eigenvalue} (with multiplicities) of the \textit{nonnormalized}  combinatorial Laplacian $L(\Gamma)$, that means, $d_{\Gamma}I_V -A(\Gamma)$. Here $A(\Gamma)$ denotes the adjacency matrix of $\Gamma$. Namely, the eigenvalues of $L(\Gamma)$ is $\lambda_1=0 \leq \lambda_2 \leq \cdots \leq \lambda_{|V|}$.
  \item In items $(1)$--$(3)$, we also define the \textit{normalized} versions $h'_n(\Gamma),\iota'_n(\Gamma), \lambda'_n(\Gamma)$ by dividing each quantity by $d_{\Gamma}$. Namely, $h'_n(\Gamma)=h_n(\Gamma)/d_{\Gamma}$, $\iota'_n(\Gamma)=\iota_n(\Gamma)/d_{\Gamma}$, and $\lambda'_n(\Gamma)=\lambda_n(\Gamma)/d_{\Gamma}$.
\end{enumerate}
Note that in the standard literature $\lambda_2$ here is written as $\lambda_1$.

\end{definition}

The parameters $h_2$, $\iota_2$, $\lambda_2$ are fundamental in spectral graph theory. They are non-zero if and only if $\Gamma$ is connected, and $2h_2/d_{\Gamma} \leq  \iota_2 \leq 2h_2$. Deeper relationships are \textit{Cheeger inequalities}, which state as follows:
\begin{itemize}
 \item (Alon--V. Milman \cite{AV}): $ \lambda_2/2 \leq h_2\leq \sqrt{2d_{\Gamma}}\sqrt{\lambda_2}$;
 \item (Bobkov--Houdr\'{e}--Tetali \cite{BHT}): $\lambda_2 \geq (\sqrt{\iota_2+1}-1)^2/4$.
\end{itemize}
The first one implies that $\lambda_2'/2 \leq h_2'\leq \sqrt{2\lambda_2'}$. However, it is well known that it is impossible to exclude contributions of $d_{\Gamma}$ from the right-hand side of the first inequality. The second inequality shows that we can bound $\lambda_2$ from below by $\iota_2$ without any dependence on $d_{\Gamma}$.

We say that an infinite sequence $\{\Gamma_m=(V_m,E_m)\}_{m\in \mathbb{N}}$ is a sequence of \textit{expanders} if the following three conditions are fulfilled:
\begin{enumerate}
  \item[$(\mathrm{ex}\ i)$] The $\sup_m d_{\Gamma_m}<\infty$; 
  \item[$(\mathrm{ex}\ ii)$] The $\lim_{m\to \infty } |V_m|=\infty$; 
  \item[$(\mathrm{ex}\ iii_2)$] The $\inf_{m} h_2(\Gamma_m)>0$. 
\end{enumerate}
The terminology ``expanders" originates from condition $(\mathrm{ex}\ iii_2)$. By Cheeger inequalities above, under condition $(\mathrm{ex}\ i)$, condition $(\mathrm{ex}\ iii_2)$ is equivalent to saying that $\inf_{m} \iota_2(\Gamma_m)>0$, as well as to saying that $\inf_{m} \lambda_2(\Gamma_m)>0$.

In terms of multi-way expansions, the following notion is defined. The notion of $2$-way expanders is identical to that of expanders.
\begin{definition}\label{defn:nexp}
For fixed $n$, a sequence of finite graphs $\{ \Gamma_m\}_{m\in \mathbb{N}}$ is called a sequence of $n$\textit{-way expanders} if conditions $(\mathrm{ex}\ i)$ and $(\mathrm{ex}\ ii)$ above; and the following condition $(\mathrm{ex}\ iii_n)$ are satisfied:
\[
(\mathrm{ex}\ iii_n):\quad \inf_{m} h_n(\Gamma_m)>0.
\]
\end{definition}

We note that under condition $(\mathrm{ex}\ i)$, condition $(\mathrm{ex}\ iii_n)$ is equivalent to saying that $\inf_{m} \iota_n(\Gamma_m)>0$, as well as to saying that $\inf_{m} \lambda_n(\Gamma_m)>0$. Indeed, this follows from $2h_n/d_{\Gamma} \leq  \iota_n \leq 2h_n$ and the following \textit{higher-order Cheeger inequality} by Lee--Oveis Gharan--Trevisan \cite[Theorem~1.1]{LGT}: 
\[
\frac{1}{2}\lambda'_n(\Gamma) \leq \rho_{\Gamma}(n) \leq O(n^2)\sqrt{\lambda'_n(\Gamma)},
\]
which is proved in providing certain spectral partitioning of graphs. 
Here $\rho_{\Gamma}(n)$ denotes the quantity $\min_{1\leq i \leq n} \max |\partial S_i|/(d_{\Gamma}|S_i|)$, where the minimum is over all collections of $n$ non-empty, disjoint subsets $S_1 ,\ldots ,S_n \subset V$. The difference between $\rho_{\Gamma}(n)$ and $h'_n(\Gamma)$ is that in the former quantity, we do not impose that $S_1,\ldots ,S_n$ is a partition of $V$ (in other words, $\bigsqcup_{i=1}^n S_i$ may be a proper subset of $V$). By the proof of \cite[Theorem~3.8]{LGT}, it is easy to see that $\rho_{\Gamma}(n) \leq h'_n(\Gamma) \leq n \rho_{\Gamma}(n)$. Therefore, their inequality above, in particular, implies that
\[
\frac{1}{2}\lambda_n(\Gamma) \leq h_n(\Gamma) \leq O(n^3)\sqrt{d_{\Gamma}} \sqrt{\lambda_n(\Gamma)}.
\]

There is a similar result independently proved by Louis--Raghavendra--Tetali--Vempala \cite{LRTV12}. There are, in addition, several related papers on spectral partitioning, see for instance \cite{KLLG} and \cite{LRTV11}. It is also worth noting that Miclo \cite{Miclo} deduced from the higher-order Cheeger inequality above the positive solution to the spectral gap conjecture for hyperbounded Markov operators.

By Lemma~\ref{lemma:gnCheeger} in Section~\ref{Section:gnCheeger}, which is a direct application of the arguments in \cite{LGT} and \cite{BHT} to the multi-way symmetric vertex isoperimetric constants,  together with the higher-order Cheeger inequality above, we have that for a fixed $n$, and for $\{\Gamma_m\}_{m\in\mathbb{N}}$, 
\[
\inf_{m}\iota_n (\Gamma_m)>0 \quad \Rightarrow \quad \inf_{m}\lambda_n(\Gamma_m)>0 \quad \Rightarrow \quad \inf_{m}h_n(\Gamma_m)>0.
\]
These three conditions are all equivalent if $\mathrm{sup}_m d_{\Gamma_m}<\infty$. However, in general case, no two of these three are equivalent. Indeed, this is not difficult to see: for instance, ``fatten" cycle graphs in appropriate ways (more precisely, consider $G_{m,k}=\mathbb{Z}/m\mathbb{Z}\times \mathbb{Z}/k\mathbb{Z}$, $S_{m,k}=\{1\}\times (\mathbb{Z}/k\mathbb{Z}\setminus \{0\})$, and $\Gamma_{m,k}=\mathrm{Cay}(G_{m,k},S_{m,k})$ for appropriately-chosen pairs $(m,k)$).

Note that $h_n, \iota_n, \lambda_n$ are non-decreasing for $n$ (for first two, observe that $|\partial (A\sqcup B)|\leq |\partial A|+|\partial B|$ and $|\delta (A\sqcup B)| \leq |\delta A|+|\delta B|$ for disjoint $A,B\subseteq V$), and hence that being $(n+1)$-way expanders are weaker than being $n$-expanders in general. This is strictly weaker.  Indeed, pick some sequence of expanders $\{\Lambda_k\}_{k\in\mathbb{N}}$ and construct a new family of graphs $\{\Gamma_m\}_{m\in \mathbb{N}}$ as follows: connect components of the disjoint union $\bigsqcup_{i=1}^{n} \Lambda_{m+i}$ each other by small number of edges (it can be done in such a way that resulting graphs are regular) and set it as $\Gamma_m$. Then $\{\Gamma_m\}_{m\in \mathbb{N}}$ are $(n+1)$-way expanders but not $n$-way expanders. Conversely, M. Tanaka \cite[Theorem~2]{Tanaka} has showed that if $h_{n+1}(\Gamma)>3^{n+1}h_{n}(\Gamma)$, then there exists a  partition of $V(\Gamma)$ into non-empty $n$ disjoint subsets $V_1,\ldots ,V_n$ such that $(h_n(\Gamma)\leq )\max_i |\partial V_i|/|V_i| \leq 3^nh_n(\Gamma)$; and that for any $i$, the induced subgraph on $V_i$ by $\Gamma$ has a ($2$-)isoperimetric constant at least $h_{n+1}(\Gamma)/3^{n+1}$. Therefore, if $h_{n+1}(\Gamma)$ is sufficiently larger than $h_{n}(\Gamma)$, then $\Gamma$ is constructed in a way similar to one above.

However, resulting graphs from the construction above do not seem homogeneous. In this point of view, Koji Fujiwara asked the following question.

\begin{question}$($K. Fujiwara$)$\label{question:F}

 For a sequence of finite connected $\mathrm{Cayley}$ graphs, does the property of being $n$-way expanders in fact imply that of being expanders for every $n$?
\end{question}

We may ask stronger question as follows:

\begin{question}\label{question:M}
\begin{enumerate}[$(1)$]
  \item   Does there exist a universal constant $C=C(n)$, depending only on $n$, such that for any finite connected Cayley graph $\Gamma$, $h_{n+1}(\Gamma)\leq C h_n(\Gamma)$ holds true?
  \item The same question with replacing $h_n$'s with $\iota_n$'s.
\end{enumerate}
\end{question}

His original idea was to translate ``thin" part to ``thick" part by the group action and to lead a contradiction if there were some counterexample to Question~\ref{question:F}. This idea is, indeed, the first step to deal with these questions.

In this paper, we provide the satisfactory answers to all of these questions. the answer to Question~\ref{question:F} is \textit{affirmative}. Item $(1)$ of Question~\ref{question:M}, however, has the \textit{negative} answer. Surprisingly, nevertheless, we answer item $(2)$  in the \textit{affirmative}. These answers follow from the following \textit{universal inequalities} for finite connected vertex-transitive graphs (observe that $\iota_{n+1}(\Gamma)\leq 2n+1$ always holds).

\begin{mthm}[Main Theorem]\label{mthm:mt}
Let $\Gamma$ be a finite connected $\mathrm{vertex}$-$\mathrm{transitive}$ graph and $2\leq n\leq |V|-1$. Then we have that 
\[
h_{n}(\Gamma)\geq \frac{h_{n+1}(\Gamma)}{10n+h_{n+1}(\Gamma)}, \quad \textrm{and} \quad  \iota_{n}(\Gamma)\geq \frac{2\iota_{n+1}(\Gamma)}{20n+\iota_{n+1}(\Gamma)}. 
\]
In particular, $\iota_{n+1}(\Gamma)\leq (11n+1)\iota_n(\Gamma)$.
\end{mthm}

Note that from the discussions above Question~\ref{question:M}, there is totally \textit{no} bound of $h_{n+1}(\Gamma)$ from above by $h_{n}(\Gamma)$ for a \textit{general} finite connected graph.

\begin{corollary}\label{cor:multiwayexp}
Let $\{\Gamma_m\}_{m\geq \mathbb{N}}$ be a sequence of finite connected $\mathrm{vertex}$-$\mathrm{transitive}$ graphs such that $\lim_{m\to \infty}|V_m|=\infty$ $($we do $\mathrm{not}$ assume that $\sup_m d_{\Gamma_m}<\infty$$)$. Then for any $n\geq 2$, $\inf_m h_{n+1}(\Gamma_m)>0$ implies $\inf_m h_{n}(\Gamma_m)>0$; and $\inf_m \iota_{n+1}(\Gamma_m)>0$ implies $\inf_m \iota_{n}(\Gamma_m)>0$. 

In particular, if $\{\Gamma_m\}_{m\in \mathbb{N}}$ are $n$-way expanders for some $n\geq 2$, then they are, in fact, expanders.
\end{corollary}

We discuss the statement of Theorem~\ref{mthm:mt} in details. For the assertions below, see Section~\ref{section:counterexample}. There we provide some counterexamples. First, we remark that since $h_{n+1}(\Gamma)\leq d_{\Gamma}$, Theorem~\ref{mthm:mt} implies that $h_{n+1}(\Gamma)\leq (10n+d_{\Gamma})h_n(\Gamma)$. However it is \textit{impossible} to avoid the contribution of the degree form the right-hand side of this inequality. Also, we note that Theorem~\ref{mthm:mt} implies that if $h_n(\Gamma)<1-\epsilon$ for some $\epsilon>0$, then 
\[
h_{n+1}(\Gamma)< \frac{10n}{\epsilon}h_n(\Gamma) 
\]
holds true. We, however, have \textit{no hope} to obtain any nontrivial estimate of $h_{n+1}(\Gamma)$ as soon as $h_n(\Gamma)\geq 1$. We, furthermore, show that this vaule $1$ is the \textit{optimal} critical value. We also warn that if we consider $h'_n$, instaed of $h_n$, or consider $\lambda'_n$, then the corresponding assertion to Corollary~\ref{cor:multiwayexp} is \textit{no longer} true. To the best knowledge of the author, similar results to above for $\lambda_n$'s seem to be open. More precisely, Li \cite{Li} showed a universal inequality for homogeneous manifold, but a naive application of his result to a vertex-transitive graph \textit{fails} to be true. The problem here may be because the vertex-transitivity of a graph can be regarded as a weaker assumption than the homogeneity of a manifold in the corresponding setting, because there is, in general, \textit{no} homogeneity of edges for a vertex-transitive graph. 

Finally, we note that, although the proofs are different, our inequalities in Theorem~\ref{mthm:mt} may have a similar philosophy to ones in \cite[Theorem 1.1]{FS} and in \cite[Theorem~1.2 and Theorem~1.6]{Funano}: their inequalities are universal, \textit{independent of dimensions of manifolds} with nonnegative Ricci curvature; and ours are universal, \textit{independent of degrees of vertex-transitive graphs}.

We, furthermore, show that if the group action $G\curvearrowright \Gamma$ possesses certain ``homogeneity", then the answer to $(1)$ of Question~\ref{question:M} is affirmative. This condition is stated in terms of \textit{primitive} group actions (for the definition of a \textit{system of imprimitivity of size $n$}, see Definition~\ref{def:systemofimprimitivity}). More precisely, we show that gaps between $n$-way isoperimetry and $(n+1)$-way one implies the existence of a system of imprimitivity of size $n$ sufficiently close to from a (fixed) realizer of $n$-way isoperimetry.

\begin{mthm}\label{mthm:imprimitive}
Let $\Gamma$ be a finite vertex-transitive graph $($possibly disconnected$)$ and $2\leq n \leq |V|-1$. If $h_{n+1}(\Gamma)> 2(n+1)h_n(\Gamma)$, then there exists decompositions $V=V_1 \sqcup V_2 \sqcup \cdots \sqcup V_n$ and $V=A_1 \sqcup A_2\sqcup \cdots \sqcup A_n$ into $n$ non-empty sets which satisfy the following properties:

\begin{enumerate}[$(a)$]
 \item The $V=V_1 \sqcup V_2 \sqcup \cdots \sqcup V_n$ is a system of imprimitivity $($of size $n$$)$  for $\mathrm{Aut}(\Gamma)\curvearrowright V$.
 \item The $V=A_1 \sqcup A_2 \sqcup \cdots \sqcup A_n$ achieves $h_n(\Gamma)$.
 \item For any $1\leq i\leq n$, $|V_i \triangle A_i| \leq \frac{4h_n(\Gamma)}{h_{n+1}(\Gamma)}|V|$.
\end{enumerate}
In fact, we may obtain $(V_i)_i$ with items $(a)$ and $(c)$ for any given $(A_i)_i$ with item $(b)$.

In particular, for $G$ a group which acts on $\Gamma$ vertex-transitively, if there exists no system of imprimitivity of size $n$ for $G\curvearrowright V$, then $h_{n+1}(\Gamma)\leq 2(n+1)h_n(\Gamma)$ holds.

The same results hold true if we replace $h_n(\Gamma)$ and $h_{n+1}(\Gamma)$, respectively, with $\iota_n(\Gamma)$ and $\iota_{n+1}(\Gamma)$.
\end{mthm}

Theorem~\ref{mthm:imprimitive} may relate to the famous problem of M. Kac, \textit{``Can one hear the shape of a drum?"}, which asks whether we can detect shapes from spectral data. A baby case of Theorem~\ref{mthm:imprimitive} is the one where $h_n=0$ and $h_{n+1}>0$. Then $\Gamma$ has exactly $n$ connected components, and we can take the associated decomposition both for $(V_i)_i$ and for $(A_i)_i$.

\

\noindent
\textbf{Plan of the proof of Theorem~\ref{mthm:mt}.} 

Despite appearance of statements in Theorem~\ref{mthm:mt}, the proof is long and complicated (though the arguments are elementary). For instance, Theorem~\ref{mthm:imprimitive}, in fact, will be proved before Theorem~\ref{mthm:mt}, and it is needed for the proof of Theorem~\ref{mthm:mt}. More precisely, we will prove Theorem~\ref{mthm:mt} in the following $4$ steps:

\noindent
\textbf{Step~1.} Take a vertex-transitive finite graph $\Gamma$. If $h_{n+1}(\Gamma)\leq 2(n+1)h_n(\Gamma)$, then it is nothing to prove. If $h_{n+1}(\Gamma)> 2(n+1)h_n(\Gamma)$, then for any subgroup $G$ of $\mathrm{Aut}(\Gamma)$ which acts transitively on $V$, we will construct a \textit{group homomorphism} $\Phi_n$ from $G$ to a symmetric group of an appropriate degree $l$, where $1\leq l\leq n$. The precise statement is as follows:

\begin{lemma}\label{lemma:Phi}
Let $\Gamma$ be a finite vertex-transitive graph $($possibly disconnected$)$, and $|V|-1\geq n\geq 2$. Take any group $G$ that acts on $\Gamma$ vertex-transitively. Assume that $h_{n+1}(\Gamma)>2(n+1)h_n(\Gamma)$. Take any partition $V=A_1 \sqcup \cdots \sqcup A_n$ which achieves $h_n(\Gamma)$. 

Then, after permuting the indices $(1,\ldots ,n)$ if necessary, there exist an integer $l \in [n]$ and a $\mathrm{group}$ $\mathrm{homomorphism}$
\[
\Phi_n \colon G \to \mathfrak{S}_l,
\]
such that $\Phi_n$ satisfies the following: for any $i\in [l]$, 
\[
|g\cdot A_i \triangle A_{\Phi_n(g)(i)}|\leq \frac{4h_n(\Gamma)}{h_{n+1}(\Gamma)}\max_{k\in [n]}|A_k|.
\]

The corresponding assertion remains true if we replace simultaneously all of the $h_{n}(\Gamma)$ and $h_{n+1}(\Gamma)$ in the assertion above, respectively,  with $\iota_{n}(\Gamma)$ and $\iota_{n+1}(\Gamma)$.
\end{lemma}
We may assume that the induced $G$-action by $\Phi_n$ on the $l$-point set $[l]$ is transitive. Then the $l$ above is uniquely determined according to the choice of $(A_1,\ldots ,A_n)$ and $G$. This homomorphism will play a central r\^{o}le in the proof of Theorem~\ref{mthm:mt}.

\noindent
\textbf{Step~2.} The proof of Theorem~\ref{mthm:mt} is divided into two cases according to the value of $l$ in the last part of Step~1: the case where $l<n$; and one where $l=n$. In this step, we treat the easier case: the former case of ``$l<n$". In fact, we will prove that this case is \textit{impossible}. This argument leads the following result, which is a weaker form of Theorem~\ref{mthm:imprimitive}, and which may be of its own interest:
\begin{theorem}\label{thm:weakimprimitive}
Let $\Gamma$ be a finite vertex-transitive graph $($possibly disconnected$)$, and $|V|-1 \geq n\geq 2$. Take any group $G$ that acts on $\Gamma$ vertex-transitively. Assume that $G$ satisfies the following condition $(\ast_n)$:
\[
(\ast_n):\ \textrm{\textit{no} action $G$ on an $n$-point set is  transitive.}
\]
Then $h_{n+1}(\Gamma)\leq 2(n+1)h_n (\Gamma)$, and $\iota_{n+1}(\Gamma)\leq 2(n+1)\iota_n(\Gamma)$ hold true.
\end{theorem}
\noindent
(Note that condition $(\ast_n)$ is characterized by the non-existence of subgroups of index $n$.) 

The key to the proof here is that if a non-empty subset $B$ of $V$ is ``almost invariant" by the transitive $G$-action, then $|V\setminus B|$ must be small. See Lemma~\ref{lem:key} for the precise statement.

\noindent
\textbf{Step~3.} Now we deal with the harder case: $l=n$. In this step, we will find a system of impripitivity $(V_1,\ldots, V_n)$ of size $n$ for the action $G\curvearrowright V$ which is ``sufficiently close" to the initially taken partition $(A_1,\ldots, A_n)$. More precisely, here we, thus, prove Theorem~\ref{mthm:imprimitive}. To do this, we consider the orbits of the characteristic functions $\chi_{A_1},\ldots ,\chi_{A_n}$ by the action $G\curvearrowright \ell_1(V)$ induced by $G\curvearrowright V$, and take level sets in appropriate sense according to coset decompositions of $G$ by $\Phi_n\colon G\to \mathfrak{S}_n$ as in Step~1. Here $\ell_1(V)$ denotes the real-valued $\ell_1$-space on $V$.

\noindent
\textbf{Step~4.} Before this step, we have not needed the assumption of that $\Gamma$ is a connected graph. Here we employ this assumption. Rough idea is as follows: by the connectedness of $\Gamma$, there exists at least one edge which connects $V_i$ and $V_j$ with some distinct $i,j$ for the $(V_1,\ldots ,V_n)$ as in Step~3. By the vertex-transitivity, we prove that for any $i\in [n]$ and for any $v\in V_i$, there exists \textit{at least $1$ edge} which connects $v$ and a vertex $w_v$ in $V\setminus V_i$; and, moreover, that we can choose $\{(v,w_v)\}_{v\in V_i}$ in such a way that a different $v\in V_i$ gives a different $w_v\in V\setminus V_i$. This deduction is straightforward  if $\Gamma$ is a Cayley graph, but in general case of vertex-transitive graphs, we use the Hall marriage theorem. Finally, by combining that observation with item $(c)$ of Theorem~\ref{mthm:imprimitive}, we establish Theorem~\ref{mthm:mt} .

Note that the \textit{critical value} $1$, as we argued in Introduction (see the paragraph below Corollary~\ref{cor:multiwayexp}), comes from this step.

\

\noindent
\textbf{Organization of this paper.} In Section~\ref{Section:gnCheeger}, we state Lemma~\ref{lemma:gnCheeger} and briefly sketch a proof of it. Section~\ref{section:counterexample} is for counterexamples to item $(1)$ of Question~\ref{question:M}. Those counterexamples, furthermore, explain our assertions in the paragraph below Corollary~\ref{cor:multiwayexp}. From Section~\ref{section:Phi} to Section~\ref{section:mt}, we follow the steps explained above to prove Theorem~\ref{mthm:mt}: in each Section, we make one step. In Section~\ref{section:ve}, we remark that item $(1)$ of Question~\ref{question:M} is resolved in the affirmative if $\Gamma$ is vertex and edge transitive. In Section~\ref{section:que}, we state some further questions.

\section{Higher-order Cheeger inequality for symmetric vertex isoperimetries}\label{Section:gnCheeger}

We give a proof of the following result, which may be regarded as a higher order Bobkov--Houdr\'{e}--Tetali inequality; 
\begin{lemma}\label{lemma:gnCheeger}
For a finite graph $\Gamma$ and $2\leq n\leq |V|$, we have that 
\[
O(n^6) \lambda_n(\Gamma) \geq \left(\sqrt{\tilde{\iota}_n(\Gamma)+1} -1\right)^2.
\]
Here $\tilde{\iota}_n(\Gamma)$ denotes the quantity $\min_{1\leq i \leq n} \max |\delta S_i|/|S_i|$, where the minimum is over all collections of $n$ non-empty, disjoint subsets $S_1 ,\ldots ,S_n \subset V$. 

In particular, we have that
\[
O(n^6) \lambda_n(\Gamma) \geq \left(\sqrt{\frac{\iota_n(\Gamma)}{n}+1} -1\right)^2.
\]
\end{lemma}
As we mentioned in Introduction, this result is a direct application of the work in \cite{LGT} to one in \cite{BHT}. For the reader's convenience, we briefly sketch the proof.

\begin{proof}
For a non-zero $f\in \ell_2(V,\mathbb{R})$, the \textit{Rayleigh quotient} $\mathrm{Ray}_{\Gamma}(f)$ of $f$ is given by 
\[
\mathrm{Ray}_{\Gamma}(f)= \frac{\sum_{(u,v)\in E}|f(u)-f(v)|^2}{\sum_{v\in V}|f(v)|^2}.
\]
Note that we consider nonnormalized one, namely, we do not divide the right-hand side by $d_{\Gamma}$. Then the following  is easily derived from arguments in \cite[Theorem~2]{BHT} (compare with Lemma~2.2 in \cite{LGT}).
\begin{lemma}\label{lem:BHT}
For any $0\ne f\in \ell_2(V,\mathbb{R})$, there exists a subset $\emptyset \ne S\subseteq \mathrm{supp}(f)$ such that 
\[4\mathrm{Ray}_{\Gamma} (f) \geq \left(\sqrt{\frac{|\delta S|}{|S|} +1}-1\right)^2.
\] 
\end{lemma}
This lemma together with Theorem~1.5 in \cite{LGT} ends our proof. For the latter item of Lemma~\ref{lemma:gnCheeger}, note that the proof of Theorem~3.8 in \cite{LGT} implies that $\tilde{\iota}_n(\Gamma) \leq \iota_n(\Gamma)\leq n\tilde{\iota}_n(\Gamma)$. Compare with the relationship between $\rho_{\Gamma}(n)$ and $h'_n(\Gamma)$ in Introduction.
\end{proof}

\section{Counterexamples to $(1)$ in Question~\ref{question:M}}\label{section:counterexample}
In this section, we will provide a construction of certain Cayley graphs that serves as counterexamples, at the same time, to  the following:
\begin{itemize}
 \item Item $(1)$ in Question~\ref{question:M};
 \item The resulting question of item $(1)$ in Question~\ref{question:M} by replacing $h_{n+1}(\Gamma)$ and $h_n(\Gamma)$, repsectively, with $\lambda_{n+1}(\Gamma)$ and $\lambda_{n}(\Gamma)$;
 \item The resulting assertion of Corollary~\ref{cor:multiwayexp} by replacing $h_{n+1}(\Gamma)$ and $h_n(\Gamma)$, repsectively, with $h'_{n+1}(\Gamma)$ and $h'_{n}(\Gamma)$;
 \item The resulting assertion of Corollary~\ref{cor:multiwayexp} by replacing $h_{n+1}(\Gamma)$ and $h_n(\Gamma)$, repsectively, with $\lambda'_{n+1}(\Gamma)$ and $\lambda'_{n}(\Gamma)$;
 \item Hope to have a similar result to one by \cite{Li} for finite connected vertex-transitive graphs;
 \item Hope to give a non-trivial bound of $h_{n+1}(\Gamma)$ when we only know that $h_n(\Gamma)$ is at most $1$.
\end{itemize}
Note that, if we know that $h_n(\Gamma)$ is at most $1-\epsilon$ for some $\epsilon >0$, then we have a non-trivial bound for $h_{n+1}(\Gamma)$, see discussions in Introduction.

First we give a counterexample for $n=2$. Let $\Lambda=\mathrm{Cay}(H,T)$ have very big $h_2$. (For instance, set $(H,T)=(\mathbb{Z}/N\mathbb{Z}, \mathbb{Z}/N\mathbb{Z}\setminus \{0\})$ for a large $N$. Then $\Lambda$ is the complete graph $K_N$.) This implies that $|T|$ is also very big. Let $G=H\times \mathbb{Z}/2\mathbb{Z}$, and set a generating set $S=(T\times \{0\})\sqcup \{(e_{H},1)\}$. Then the Cayley graph $\Gamma=\mathrm{Cay}(G, S)$ is a counterexample (note that this graph is the graph product of $\Lambda$ and $\mathrm{Cay}(\mathbb{Z}/2\mathbb{Z}, \{1\})$). Indeed, by decomposing as $G=(H \times \{0\})\sqcup (H \times \{1\})$, we have that $h_2(\Gamma)\leq 1$. However, Lemma~1 in \cite{Tanaka} implies that $h_3(\Gamma)\geq h_2(\Lambda)$, and this shows that we can have $h_3(\Gamma)$ as large as we wish with appropriate choices of $(H,T)$ (for instance, let $N \to \infty$). 

To show that $1$ is the critical value for $h_n$ to bound $h_{n+1}$ (see Section~\ref{sec:intro}), we modify this construction if $n\geq 3$. Take a dihedral group $D_{n}=\langle a,b\mid a^2=b^2=(ab)^n=e_{D_{n}}\rangle$, and from $(H,T)$ construct $(G,S)$ as follows: $G=H\times D_{n}$, and $S=(T\times \{e_{D_{n}},a\}) \sqcup \{(e_{H},b)\}$. Then for $\Gamma=\mathrm{Cay}(G,S)$, a similar argument to above tells us that $h_{n}(\Gamma)\leq 1$; but that $h_{n+1}(\Gamma)$ can be arbitrarily big. To see these assertions, more precisely, decompose $V(\Gamma)=G$ as $G=\bigsqcup_{i=0}^{n-1} (H \times \{(ab)^i, (ab)^ia\})$. Then $h_n(\Gamma)\leq 1$, and \cite[Lemma1]{Tanaka} shows that $h_{n+1}(\Gamma)\geq h_2(\Lambda)$.

To see that these are also counterexamples to the corresponding question to $\lambda_n$'s in Question~\ref{question:M}, apply \cite[Lemma6]{Tanaka}. In particular, we cannot naively apply Li's results in \cite{Li} on $\lambda'_n$'s to the case of finite connected vertex-transitive graphs (because otherwise \cite[Theorem~11]{Li}, in particular, would imply that $\lambda'_{n+1}(\Gamma)<5 \lambda'_n(\Gamma)$). We, in addition, note that if we consider weighted cases, then the corresponding assertions in Corollary~\ref{cor:multiwayexp} \textit{fail} to be true. (For the definition of $h'_n$ and $\lambda'_n$ for weighted graphs, see \cite{LGT}.) More precisely, in that case, we put a weight on $S$. If we put very small weight on $(e_{H},b)$ relative to the other elements in $S$ in the example above, then this construction serves as counterexamples to the assertions both on weighted $h'$ and weighted $\lambda'$. These counterexamples may be constructed even in such a way of that the degrees of the graphs are uniformly bounded. 

\section{Step~1: Construction of group homomorphisms from gaps between multi-way isoperimetires}\label{section:Phi}

From this section to Section~\ref{section:mt}, we focus on the proof of Theorem~\ref{mthm:mt}. In this section, we verify Lemma~\ref{lemma:Phi}. The following lemma is obvious, and we will employ it without mentioning throughout the present paper.

\begin{lemma}\label{lem:triv}
Let $\Gamma=(V,E)$ be a finite graph and $U_1 ,U_2\subseteq V$.
\begin{enumerate}[$(1)$]
  \item Let $U_1\cap U_2=\emptyset$. Then for any $g\in \mathrm{Aut}(\Gamma)$, $|\partial (g\cdot U_1, g \cdot U_2)|=|\partial (U_1,U_2)|$. In particular, $|\partial (g\cdot U_1)|=|\partial U_1|$.
  \item For $U_1\subseteq U_1'\subseteq V$ and $U_2\subseteq U_2'\subseteq V$ with $U_1'\cap U_2'=\emptyset$, $|\partial (U_1,U_2)| \leq |\partial (U_1',U_2')|$.
  \item We have that $|\partial (U_1\cap U_2)| \leq |\partial U_1| + |\partial U_2|$. Moreover, $|\partial (U_1\cap U_2)| \leq |\partial U_1|+ |\partial (U_1\setminus U_2, U_1 \cap U_2)| $. 
  \item If $U\subseteq V$ is partitioned as $U=\bigsqcup_{j=1}^k U_j$, then $|\partial U| \leq  \sum_{j=1}^k |\partial U_j|$.
\end{enumerate}
All of the corresponding statements remain true if we replace all $\partial$ with $\delta$ in the setting above.
\end{lemma}

Before proceeding in the proof of Lemma~\ref{lemma:Phi}, we present rough idea. Take any partition $(A_1,\ldots,A_n)$ of $V$ which achieves $h_n(\Gamma)$, and any $G \leqslant \mathrm{Aut}(\Gamma)$ which acts transitively on $V$. Fix any $g\in G$ and $j\in [n]$ such that $|A_j|$ is not small among $|A_1|, \ldots , |A_n|$. Then, because of the gap between $h_{n+1}(\Gamma)$ and $h_n(\Gamma)$, for any $k\in [n]$, either $g\cdot A_j$ is ``close to" $A_k$; or $g\cdot A_j$ is ``almost disjoint from" $A_k$. From some quantitative estimate, we are able to show that, in fact, for the pair $(g,j)$ above, there exists a \textit{unique} $k=k(g,j)\in [n]$ which fulfills the first option. We, thus, obtain a map $I_j\colon G \to [n]$ which maps $g$ to $k=k(g,j)$. 

We, furthermore, sketch the way to obtain a homomorphism law in our constrution. Take $g,g'\in G$ and set $j=I_i(g')$ and $k=I_j(g)$. Then, by the definition of $I_i$ and $I_j$, $g'\cdot A_i$ is ``close to" $A_j$; and $g\cdot A_j$ is ``close to" $A_k$. From them, we wish to conclude that $(gg')\cdot A_i$ is ``close to" $A_k$ by considering the composition of the multiplications $A_i \mapsto g' \cdot A_i \mapsto (gg')\cdot A_i$. However, there is one problem in this deduction:  the error between $(gg')\cdot A_i$ and $A_k$, in general, might get bigger than the admissible error in the first option for the pair $(gg',i)$. We overcome this difficulty by the following  key observation: \textit{there are only two options}: $(gg')\cdot A_i$ and $A_k$ are ``\textit{close}"; or they are ``\textit{almost disjoint}." Even if the error between $(gg')\cdot A_i$ and $A_k$ might grow, under our assumption, it is impossible to meet the second option. Therefore, \textit{the $k$ above must satisfy the first option} for $(gg',i)$. This argument shows that $I_{i}(gg')=I_{I_{i}(g')}(g)$. From this equality, we can construct a \textit{group homomorphism} $\Phi_n \colon G\to \mathfrak{S}_l$ for an appropriate $l\in [n]$, after permuting indices $1,\ldots ,n$ if necessary.

\begin{proof}[Proof of Lemma~$\ref{lemma:Phi}$]

We will only show the assertion for $h_n$ (the proof for $\iota_n$ goes exactly along the same way). 
 
First, note that the assumption of that $h_{n+1}(\Gamma)>2(n+1)h_{n}(\Gamma)$, in particular, implies that $6h_n(\Gamma)<h_{n+1}(\Gamma)$. Let $(A_1,\ldots ,A_n)$ be a (non-empty) $n$-partition of $V$ which achieves $h_n(\Gamma)$. Without loss of generality, we may assume that $|A_1|$ is the largest among $|A_1|,\ldots,|A_n|$. 

Secondly, fix $g \in G$. For each $1\leq k\leq n$, decompose $V$ into $g^{-1} \cdot A_k \cap A_1$, $A_1 -g^{-1} \cdot A_k$, and $A_2,\ldots ,A_n$. Because 
\[
|\partial(g^{-1} \cdot A_k \cap A_1, A_1 -g^{-1} \cdot A_k)|\leq |\partial A_k|\leq h_n(\Gamma)|A_k|\leq h_n(\Gamma)|A_1|,
\]
we have that 
\[
|\partial(g^{-1} \cdot A_k \cap A_1)|\leq h_n(\Gamma)|A_1|+|\partial(A_1)|\leq 2h_n(\Gamma)|A_1|,
\] 
and that $|\partial( A_1 -g^{-1} \cdot A_k)|\leq 2h_n(\Gamma)|A_1|$.
From the condition of $h_{n+1}(\Gamma)$, we conclude the following: for fixed $g \in G$, for each $1\leq k\leq n$, either of the following $(i)_1$ and $(ii)_1$ holds true:
\begin{enumerate}
   \item[$(i)_1$]: $|g\cdot A_1\cap A_k|\geq \left(1-\frac{2h_n(\Gamma)}{h_{n+1}(\Gamma)}\right)|A_1|$;
   \item[$(ii)_1$]: $|g\cdot A_1\cap A_k|\leq \frac{2h_n(\Gamma)}{h_{n+1}(\Gamma)}|A_1|$.
\end{enumerate}
(Note that if either of two sets in the decomposition is empty, then the assertion above trivially holds.) Because $4h_n(\Gamma)<h_{n+1}(\Gamma)$, these two options are exclusive. 

Thirdly, we claim that for each $g\in G$, there \textit{exists} a \textit{unique} $k\in [n]$ which satisfies $(i)_1$. Indeed, if there exist at least $2$ such $k$'s, then 
\[
|A_1|=\left|\bigsqcup_{k=1}^n( g\cdot A_1\cap A_k)\right| \geq 2\left(1-\frac{2h_n(\Gamma)}{h_{n+1}(\Gamma)}\right)|A_1|,
\] 
but it is absurd. Also if there is no such $k$, then all $k$ satisfies $(ii)_1$ and hence 
\[ 
|A_1|=\left|\bigsqcup_{k=1}^n (g\cdot A_1\cap A_k)\right| \leq 2n\frac{h_n(\Gamma)}{h_{n+1}(\Gamma)}|A_1| <|A_1|,
\] 
and it is again a contradiction. Thus, we can define a map which send each $g\in G$ to the unique index $k=k(g)$ for which $(i)_1$ is satisfied, and set this map as $I_1\colon G \to [n]$.

By changing the indices $2,\ldots,n$ if necessary, we may assume that there exists $l\in [n]$ such that $\mathrm{Im} (I_1)=[l]$ (note that $I_1(e)=1$). An important observation is that for any $2\leq j\leq l$, we have that
\[
|A_j| \geq \left( 1-\frac{2h_n(\Gamma)}{h_{n+1}(\Gamma)}\right) |A_1|\left(\geq \frac{n}{n+1}|A_1|\right)\textrm{\ \ $\cdots$} \textrm{$(\diamond)$}
\]
because $I_1^{-1}(j)\ne \emptyset$. In the next paragraph, we proceed to an argument which is needed if $l\geq 2$. If $l=1$, then we do not do anything there.

Fix $2\leq j\leq l$. For fixed $g\in G$, in a similar argument to one above, we have that for any $1\leq k\leq n$, 
\begin{align*}
|\partial(g^{-1} \cdot A_k \cap A_j)|\leq h_n(\Gamma)(|A_1|+|A_j|),\ |\partial( A_j -g^{-1} \cdot A_k)|\leq h_n(\Gamma)(|A_1|+|A_j|).
\end{align*}
Hence, we similarly conclude that (for each $g\in G$ and) for any $1\leq k\leq n$, either of the following $(i)_j$ and $(ii)_j$ holds true:
\begin{enumerate}
   \item [$(i)_j$]: $|g\cdot A_j\cap A_k|\geq  |A_j|- \frac{h_n(\Gamma)}{h_{n+1}(\Gamma)}(|A_1|+|A_j|)\left( \geq |A_j|-\frac{2h_n(\Gamma)}{h_{n+1}(\Gamma)} |A_1|\right)$;
   \item [$(ii)_j$]: $|g\cdot A_j\cap A_k|\leq \frac{h_n(\Gamma)}{h_{n+1}(\Gamma)}(|A_1|+|A_j|)\left( \leq \frac{2h_n(\Gamma)}{h_{n+1}(\Gamma)} |A_1|\right)$.
\end{enumerate}
Note that from $(\diamond)$ these two options are exclusive. In a similar argument to the one above, we can show that (for a fixed $2\leq j\leq l$ and) for each $g \in G$, there exists a \textit{unique} $k$ which satisfies $(i)_j$. Thus for each $2\leq j\leq l$, we get a map $I_j\colon G \to [n]$ by sending $g\in G$ to $k$ for which $(i)_j$ is satisfied. We will show the following lemma:

\begin{lemma}\label{lem:Sl}
Let $1\leq j\leq l$.
\begin{enumerate}[$(1)$]
   \item The $\mathrm{Im}I_j$ satisfies that  $\mathrm{Im}I_j \subseteq [l]$.
   \item For each $g \in G$, we define $\sigma_{g}\colon [l]\to [l]$ by $\sigma_{g}(j)=I_j(g)$. Then for any $g \in G$, $\sigma_{g}\in \mathrm{Aut}([l])$$\cong \mathfrak{S}_l$.
    \item For any $g,g'\in G$, $\sigma_{g}\sigma_{g'}=\sigma_{gg'}$.
\item If $I_j(g)=k$, then we have that $|A_k\triangle g \cdot A_j| \leq \frac{h_n(\Gamma)}{h_{n+1}(\Gamma)}(2|A_1|+|A_j|+|A_k|)\left( \leq \frac{4h_n(\Gamma)}{h_{n+1}(\Gamma)}|A_1|\right)$.
\end{enumerate}
\end{lemma}

\begin{proof}(Lemma~\ref{lem:Sl})

\begin{enumerate}[$(1)$]
  \item Suppose, to the contrary, that there exists $k>l$ such that $k\in \mathrm{Im}I_j $. Because $I_j^{-1}(k)\ne \emptyset$, there exists $g\in G$ such that \[
|g\cdot A_j - A_k| \leq  \frac{2h_n(\Gamma)}{h_{n+1}(\Gamma)}|A_1|. 
\]
Because $j\in \mathrm{Im}I_1$, there again exists $g'\in G$ such that 
\[
|g g'\cdot A_1 -g\cdot A_j |= |g'\cdot A_1 -A_j |\leq  \frac{2h_n(\Gamma)}{h_{n+1}(\Gamma)}|A_1|.
\] 
By combining these two inequalities, we obtain that $|gg' \cdot A_1 -A_k| \leq  \frac{4h_n(\Gamma)}{h_{n+1}(\Gamma)}|A_1|$. Recall that, by assumption, in particular $6h_n(\Gamma)<h_{n+1}(\Gamma)$ holds. This implies that $k$ cannot satisfy option $(ii)_1$ for $gg'$: otherwise we must have that $|A_1|<|A_1|$. Therefore, $I_1(gg')=k$ but this is a contradiction.
 \item In a similar argument to one in the proof of $(1)$, we have that for any $g \in G$,
$\sigma_{g}\sigma_{g^{-1}}=\sigma_{g^{-1}}\sigma_{g}=\mathrm{id}_{\{ 1,\ldots ,l\}}$. Hence $\sigma_{g}\in \mathrm{Aut}([l])$.
 \item This can be also showed in  a similar argument to one in the proof of $(1)$.
 \item First, because $I_j(g)=k$, we have that 
\[
|g\cdot A_j -A_k| \leq \frac{h_n(\Gamma)}{h_{n+1}(\Gamma)}(|A_1|+|A_j|). 
\]
Secondly, from item $(2)$ in this lemma, we have that $I_k(g^{-1})=j$ and hence that 
\[
| A_k -g\cdot A_j  | =| g^{-1}\cdot A_k -A_j  |\leq \frac{h_n(\Gamma)}{h_{n+1}(\Gamma)}(|A_1|+|A_k|). 
\]
By combining these two inequalities, we  get the conclusion.
\end{enumerate}
\end{proof}
Finally, define the desired \textit{group homomorphism} $\Phi_n$ by
\[
\Phi_n\colon G \to \mathfrak{S}_l;\quad g\mapsto \sigma_{g},
\]
that ends our proof of Lemma~\ref{lemma:Phi}.
\end{proof}

\section{Step~2: exculsion of the case where $l<n$}\label{section:npoints}

In this section, we keep the assumptions on $n$, $\Gamma$, $(A_1,\ldots, A_n)$, and $G$ in Lemma~\ref{lemma:Phi}. Assume that $h_{n+1}(\Gamma)>2(n+1)h_n(\Gamma)$. (The case where $\iota_{n+1}(\Gamma)>2(n+1)\iota_n(\Gamma)$ can be treated in a similar way, and we omit it.) Then, by Lemma~\ref{lemma:Phi}, after permuting the indices $1,\ldots,n$, there exists an $l\in [n]$ such that we have the group homomorphism
\[
\Phi_n\colon G\to \mathfrak{S}_{l},
\]
with respect to which $G$ acts on $[l]$ transitively. Here we may assume that $|A_1|=\max_{i\in [n]}|A_i|$, and then inequality $(\diamond)$ in Section~\ref{section:Phi} holds true.

As we argued in Introduction, the goal in this section is to show that, under this assumption, $l$ must coincide with $n$. In the proof of it, we use the following  key lamma.

\begin{lemma}\label{lem:key}
Let $\varepsilon>0$. Let a finite group $H$ act on a finite set $W$ transitively. Assume that a $\mathrm{non}$-$\mathrm{empty}$ subset $C\subseteq W$ satisfies that for any $h\in H$, $| C\triangle h \cdot C|\leq \varepsilon |C|$. Then, we have that $|W\setminus C|\leq \frac{\varepsilon}{2}|W|$. In particular, if $|C|\leq |W|/2$, then $\varepsilon\geq 1$.
\end{lemma}

This lemma may be showed in a purely combinatorial argument. Instead, we here give a functional analytic proof, because we will use the same spirit of ``dualizing" (namely, to consider the characteristic function instead of a set itself) in Section~\ref{section:imprimitive}.

\begin{proof}(Lemma~\ref{lem:key})
On the (finite dimensional) Banach space $\ell_{1,0}(H)=\{\xi \in \ell_1(W):\sum_{w\in W}\xi(w)=0\}$ with the $\ell_1$-norm, a linear isometric $H$-representation $\pi$ is induced by the permutations $H\curvearrowright W$. Namely, we set as 
$\pi(h)\xi(x)=\xi(h^{-1}\cdot x)$. Note that there does not exist a nonzero $\pi(H)$-invariant vector in $\ell_{1,0}(W)$ because $H\curvearrowright W$ is transitive. Set $\xi=|W\setminus C|\chi_C- |C|\chi_{W\setminus C}$$(=|W|\chi_C-|C|\mathbf{1})\in \ell_{1,0}(W)$, where $\chi_A$ denotes the characteristic function of $A$ and $\mathbf{1}$ means the constant $1$ function. Then $\|\xi\|=2|C||W\setminus C|$, where $\|\cdot\|$ is the $\ell_1$-norm. By the assumption of the lemma, for any $h\in H$, $\|\xi-\pi(h)\xi\|\leq \varepsilon |C||W|$.

Set $\eta=|H|^{-1}\sum_{h\in H}\pi(h)\xi \in \ell_{1,0}(W)$. Because $\eta$ is $\pi(H)$-invariant, $\eta$ must be $0$. We also have that
\[
\|\xi -\eta\|=\frac{1}{|H|}\left\|\sum_{h\in H}(\xi-\pi(h)\xi)\right\|\leq \frac{1}{|H|}\sum_{h\in H}\|\xi-\pi(h)\xi\|.
\]
Therefore, we conclude that $2|C||W\setminus C|\leq \varepsilon |C||W|$.
\end{proof}

Our proof of the goal above is by the way of contradiction. More precisely, we assume that $l<n$. Then, $B=V\setminus \bigsqcup_{i=1}^l A_l$ is \textit{not} empty. We will show that $|B|<|V|/3$, and that for any $g\in G$, $|g\cdot B \triangle B| <|B|$ holds true. Then, Lemma~\ref{lem:key} derives the desired contradiction.

\begin{proof}[Proof of the equality $l=n$]
We stick to the setting at the beginning of this section. Suppose, to the contrary, that $l<n$. We set $A, B\subseteq V$ as $A=\bigsqcup_{j=1}^lA_j$, and $B=V\setminus A$; and rename $A_{l+1},\ldots, A_{n}$, respectively, as $B_{1},\ldots, B_{n-l}$. Note that  $B$ is \textit{non-empty}. We also note that from $(\diamond)$, $|A|=|A_1|+\sum_{j=2}^l|A_j|\geq \frac{ln}{n+1}|A_1|$.

First, we claim that $|B|<2|A|$. Indeed, For each $1\leq j\leq l$, we have that for each $g \in G$, $|g\cdot A_j \cap B|\leq \frac{2h_n(\Gamma)}{h_{n+1}(\Gamma)}|A_1|$ (consider $\Phi_n(g)(j)$), and that 
\[
|g\cdot A \cap B|\leq \frac{2lh_n(\Gamma)}{h_{n+1}(\Gamma)}|A_1|\leq \frac{2l(2n+3)}{l(2n+1)}\frac{h_n(\Gamma)}{h_{n+1}(\Gamma)}|A| <\frac{1}{3}|A|. 
\]
Hence, we obtain that for any $g \in G$, $|g\cdot A \triangle A|< \frac{2}{3}|A|$. 
By Lemma~\ref{lem:key}, we conclude that $|B|=|V\setminus A|< \frac{1}{3}|V|$.

In what follows, we will show our second claim: for any $g \in G$, $|B\triangle g\cdot B| <|B|$. To see this, fix $g\in G$. For any $1\leq k\leq n-l$, we have that 
\[
|\partial(g\cdot  B_k \cap A, A- g\cdot  B_k)|\leq h_n(\Gamma)|B_k|
\] 
and that 
\[|\partial(g\cdot  B \cap A)|\leq \sum_{k=1}^{n-l} h_n(\Gamma)|B_k|+\sum_{m=1}^{n-l}|\partial B_m|\leq 2h_n(\Gamma)|B|.
\]
(Here we apply item $(3)$ of Lemma~\ref{lem:triv} to the case where $U_1=A$ and $U_2= g\cdot B$, and apply item $(4)$ of Lemma~\ref{lem:triv} to make estimate of $|\partial A|=|\partial B|$.) 
Hence for any $1\leq j\leq l$, 
\begin{align*}
|\partial(A_j-g \cdot B)|&\leq 2h_n(\Gamma)|B|+ h_n(\Gamma)|A_j|< h_n(\Gamma)|A|+ h_n(\Gamma)|A_j| \\
& \leq (l+1)h_n(\Gamma)|A_1| \leq nh_n(\Gamma)|A_1|
\end{align*}
(recall that we have verified that $|A|>2|B|$). We also observe that, according to $g$ and $j$, $j'=\Phi_n({g}^{-1})(j)$ satisfies that $|g \cdot A_{j'}\cap A_j| \geq |A_j|-\frac{2h_n(\Gamma)}{h_{n+1}(\Gamma)} |A_1|$. This implies that
\begin{align*}
| A_j-g \cdot B| \geq |A_j|-\frac{2h_n(\Gamma)}{h_{n+1}(\Gamma)} |A_1|  \geq \left(\frac{n}{n+1}-\frac{2h_n(\Gamma)}{h_{n+1}(\Gamma)} \right)|A_1|.
\end{align*}
Therefore, we have the following inequalities:
\begin{align*}
\frac{|\partial(A_j-g \cdot B)|}{| A_j-g \cdot B|}< \frac{nh_n(\Gamma)}{\frac{n}{n+1}-\frac{2h_n(\Gamma)}{h_{n+1}(\Gamma)}} < (2n+2)h_n(\Gamma) < h_{n+1}(\Gamma).
\end{align*}
Finally, for $g\in G$, we decompose $V$ into $(n+1)$ disjoint subsets: $g\cdot B \cap A$; $B_1;\ldots ;B_{n-l}$; and  $A_j-g\cdot B$ ($j\in [l]$). Note that the argument above shows that $A_j-g\cdot B\ne\emptyset$ for all $j$. If $g\cdot B \cap A=\emptyset$, then $| B\triangle g \cdot B|=0$ and we are done. Hence, we may assume that all of the $(n+1)$ subsets are non-empty. Then from the condition of $h_{n+1}$, at least one subset $C$ of these $(n+1)$ subsets must satisfy that $\frac{|\partial C|}{|C|} \geq h_{n+1}(\Gamma)$. 
However by construction, neither of $B_1,\ldots ,B_{n-l}$ satisfies this condition. From the inequalities above, all of the $A_j-g\cdot B$'s, $1\leq j\leq l$ also fail to do so. Therefore, $C=g\cdot B \cap A$ must satisfy that condition. This implies that $|g \cdot B \cap A| \leq \frac{2h_n(\Gamma)}{h_{n+1}(\Gamma)}|B|$, and hence, we have that 
\[
|g\cdot B \triangle B| \leq \frac{4h_n(\Gamma)}{h_{n+1}(\Gamma)}|B|<|B|.
\]
This completes the proof of our second claim. These two claims contradict Lemma~\ref{lem:key} because $B$ is non-empty. Therefore, $l$ must equal $n$.
\end{proof}

\begin{remark}
By taking the contraposition of the statement ``$l=n$", we obtain Theorem~\ref{thm:weakimprimitive}. Theorem~\ref{thm:weakimprimitive} implies that, if $G$ is subject to condition $(\ast_n)$, then item $(1)$ of Question~\ref{question:M} resolves in the affirmative, even if $\Gamma$ is disconnected.

Because we may take any $G\leqslant \mathrm{Aut}(\Gamma)$ arbitrarily as long as $G$ acts on $V$ transitively, Theorem~\ref{thm:weakimprimitive} applies to several Cayley graphs for certain $n$. One example is a Cayley graph of  $\mathfrak{S}_N$ for $N\geq 5$. Because $\mathfrak{S}_N$ has only three normal subgroups: $\{e\}$; the alternating group $\mathfrak{A}_N$ of degree $N$; and $\mathfrak{S}_N$ itself, Theorem~\ref{thm:weakimprimitive} applies to \textit{any} Cayley graph of $\mathfrak{S}_N$ for all $3\leq n \leq N-1$. 
\end{remark}

\section{Step~3: Proof of Theorem~\ref{mthm:imprimitive}}\label{section:imprimitive}

We recall the definition of a \textit{system of imprimitivity} \textit{of size }$n$.

\begin{definition}\label{def:systemofimprimitivity}
Let $G\curvearrowright V$ be a finite group action on a finite set that is transitive. Let $n\geq 2$. A non-empty decomposition $(V_1,\ldots ,V_n)$ of $V$ ($V=V_1\sqcup \cdots \sqcup V_n$) is called a \textit{system of imprimitivity} (\textit{of size} $n$) if for any $g \in G$ there exists $\sigma_{g}\in  \mathfrak{S}_n$ such that $g \cdot V_i=V_{\sigma_{g}(i)}$ for all $1\leq i\leq n$. Each $V_i$ is called a \textit{block}.
\end{definition}

Intuitively, if a system of imprimitivity exists, then the group action does not ``break" the partitions given by blocks. It is well-known that $G\curvearrowright V$ admits a system of imprimitivity of size $n$ if and only if there exists a subgroup of $G$ of index $n$ between $G$ and a point stabilizer. For instance, compare with  \cite[Theorem~1.5A]{BookDM}.

We stick to the setting in the first paragraph in Section~\ref{section:npoints}. We choose $G=\mathrm{Aut}(\Gamma)$. (Again, we only discuss the case on $h_n$.)  Here, we are intended to prove Theorem~\ref{mthm:imprimitive}. The strategy is as follows:  In Section~\ref{section:npoints}, we have proved that $l=n$, and hence 
\[
\Phi_n \colon G \to \mathfrak{S}_n
\]
induces the transitive $G$-action on $[n]$. From this, for an each fixed $j\in [n]$, $G$ can be decomposed, in a ``coset decomposition", into $n$ subsets according to the image $\Phi_n(g)(j)$ of $j$. Then, we employ the ``dual" picture (see the paragraph below Lemma~\ref{lem:key}). We translate the characteristic functions $\chi_{A_1},\ldots ,\chi_{A_n}$, respectively, by members in appropriate cosets of $G$; and for each $i\in [n]$, take the average of the translates of $\chi_{A_i}$. Finally, by taking suitable level sets of such averaged functions, we obtain a system of imprimitivity of size $n$ for $G\curvearrowright V$ that is sufficiently ``close to" the original partition $(A_1,\ldots ,A_n)$.

\begin{proof}[Proof of Theorem~$\ref{mthm:imprimitive}$]

Keep the setting in the paragraph above. For each $(i,j)\in [n]\times [n]$, we define $G_{i,j}$ as $G_{i,j}=\{g\in G: \Phi_n(g)(j)=i \}$ (the condition on $G_{i,j}$ may be understood as ``$i \stackrel{\Phi_n(g)}{\longmapsfrom} j$"). Note that $|G_{i,j}|= |G|/n$. 

Consider the Banach space $\ell_1(V)$ with the $\ell_1$-norm, and denote by $\rho$ the isometric linear representation of $G$ on $\ell_1(V)$ by permutations. More precisely. $\rho(g)\eta(v)=\eta(g^{-1}\cdot v)$. For each $(i,j)\in [n]\times [n]$, define $M_{i,j}$ as the averaging operator on $\rho (G_{i,j})$, namely, 
\[
M_{i,j}\eta= \frac{1}{|G_{i,j}|}\sum_{g \in G_{i,j}}\rho(g)\eta(=\frac{n}{|G|}\sum_{g \in G_{i,j}}\rho(g)\eta). 
\]
Note that for any $i,j,k\in [n]$ and for any $g\in G_{i,j}$, $\rho(g)M_{j,k}=M_{i,k}$ holds. 

Set $\xi_1=\chi_{A_1},\ldots ,\xi_n=\chi_{A_n}$, and for each $i\in [n]$, define
\[
\zeta_i= \frac{1}{n}(M_{i,1}\xi_1+ M_{i,2}\xi_2+\cdots +M_{i,n}\xi_n).
\]
We claim the following:
\begin{enumerate}[$(1)$]
  \item The $\sum_{i=1}^n \zeta_i=\mathbf{1}$ and $\zeta_i(v)\in [0,1]$ for any $v\in V$ and $i\in [n]$.
  \item For any $g \in G_{i,j}$, $\rho(g) \zeta_j=\zeta_i$.
  \item For any $i$, $\|\zeta_i -\xi_i\|\leq \frac{n-1}{n}\frac{4h_n(\Gamma)}{h_{n+1}(\Gamma)}|A_1|\leq \frac{n^2-1}{n(n^2+1)}\frac{4h_n(\Gamma)}{h_{n+1}(\Gamma)}|V|$. Here $\|\cdot \|$ means the $\ell_1$-norm.
\end{enumerate}
Indeed, item $(1)$ follows from $\sum_{i=1}^n\xi_i=\mathbf{1}$ and the construction. Item $(2)$ is by $\rho(g)M_{j,k}=M_{i,k}$ and $|G_{i,j}|=|G|/n$. Item $(3)$ can be confirmed by Lemma~\ref{lemma:Phi}, the triangle inequality, and inequality $(\diamond)$ in Section~\ref{section:Phi}.

Finally, define $V_1,\ldots ,V_n$ by setting for every $i\in [n]$
\[
V_i=\left\{v\in V: \zeta_i(v)>\frac{1}{2}\right\}.
\]
We will show that $(V_1,\ldots ,V_n)$ and $(A_1,\ldots, A_n)$ satisfy all of the conclusions $(a)$--$(c)$ in Theorem~\ref{mthm:imprimitive}. First, we discuss $(a)$ and $(b)$. Item $(b)$ is by definition. To see $(a)$, observe that $V_1\ne \emptyset$ by items $(1)$ and $(3)$ above, and that for any $g\in G_{i,j}$, $g\cdot V_j=V_i$ by item $(2)$. Also, $V_i$'s are pairwise disjoint because otherwise $\sum_{i=1}^n \zeta_i\ne \mathbf{1}$. By the transitivity of the action, we see that $\bigcup_{i=1}^n V_i=V$. Hence $(V_1,\ldots, V_n)$ is a decomposition of $V$, and moreover is a system of imprimitivity of size $n$. 

Finally, we deal with the proof of item $(c)$. Because $\zeta_i$ is $\rho(G_{i,i})$-invariant (by item $(2)$), items $(1)$ and $(3)$ shows that for every $i\in [n]$ and $v\in V$,
\[
\zeta_i(v) \in \left[0, \frac{n^2-1}{(n^2+1)}\frac{4h_n(\Gamma)}{h_{n+1}(\Gamma)}\right]\cup \left[1-\frac{n^2-1}{(n^2+1)}\frac{4h_n(\Gamma)}{h_{n+1}(\Gamma)},1\right] \quad (\subseteq \mathbb{R})
\]
holds (note that $\chi_{A_i}$ takes values only in $\{0,1\}$). Therefore for every $i\in [n]$,
\begin{align*}
|V_i \triangle A_i| \leq \frac{\frac{n^2-1}{(n^2+1)}\frac{4h_n(\Gamma)}{h_{n+1}(\Gamma)}}{1-\frac{n^2-1}{(n^2+1)}\frac{4h_n(\Gamma)}{h_{n+1}(\Gamma)}}|V|  \leq \frac{4h_n(\Gamma)}{h_{n+1}(\Gamma)}|V|,
\end{align*}
as desired.
\end{proof}

\section{Step~4: end game}\label{section:mt}

We are now in position to bring the proof of Theorem~$\ref{mthm:mt}$ to an end. We, finally, make use of the assumption of the connectivity of the $\Gamma$. For idea, see Introduction.

\begin{proof}[Proof of Theorem~$\ref{mthm:mt}$]
First we prove the inequality for $h_n$'s. If $h_{n+1}(\Gamma)\leq 2(n+1)h_{n}(\Gamma)$, then we are done. Otherwise, by Theorem~\ref{mthm:imprimitive}, we may take $V=V_1\sqcup \cdots \sqcup V_n$ and $V=A_1\sqcup \cdots \sqcup A_n$ in the statement. For $G=\mathrm{Aut}(\Gamma)$, take $G_{i,j}$ for $(i,j)\in [n]\times [n]$ in Section~\ref{section:imprimitive}. 

Now we employ the assumption of that $\Gamma$ is \textit{connected}. This implies that for any $i$, there exist $v_i\in V_i$ and an edge which connects $v_i$ to a vertex $w_i$ lying in other $V_j$. Then, by translating by $G_{i,i}$-action, we observe that any $v\in V_i$, there exists at least one edge $(v,w)$ with $w\in V\setminus V_i$. 

Here we claim that we can take $w=w_v$ in such a way that a different $v\in V_i$ gives a different $w$. This claim is trivial if $\Gamma$ is a Cayley graph: switch $G$ from $\mathrm{Aut}(\Gamma)$ to the original group. To prove the claim in a general vertex-transitive graph, we will apply the Hall marriage theorem, as follows. Take a pair $v_i\in V_i$ and $w_i\in V\setminus V_i$ as above and fix them. Take any $\emptyset \neq K \subseteq V_i$ and define $H_{K}=\{g\in G: g\cdot v_i \in K\}\subseteq G_{i,i}$ (we may replace $G$ with $G_{i,i}$ above). Then from the construction, we have that 
\[
|K|= \sum_{v\in K}\frac{|\{g\in G_{i,i}:g\cdot v_i=v\}|}{|\mathrm{Stab}_{v}\cap G_{i,i}|}.
\]
Here for $y\in V$, $\mathrm{Stab}_y\leq G$ denotes the stabilizer of $y$ for $G\curvearrowright V$. Because $v\in V_i$ and the $G$-action is transitive, we have that $\mathrm{Stab}_{x}\leq G_{i,i}$ for any $x\in V_i$ and that $|\mathrm{Stab}_{y}|=|\mathrm{Stab}_{v_i}|$ for any $y\in V$. We obtain that
\[
|K|=\sum_{v\in K}\frac{|\{g\in G_{i,i}:g\cdot v_i=v\}|}{|\mathrm{Stab}_{v_i}|}=\frac{\sum_{v\in K}|\{g\in G_{i,i}:g\cdot v_i=v\}|}{|\mathrm{Stab}_{v_i}|}=\frac{|H_K|}{|\mathrm{Stab}_{v_i}|}.
\]
Let $V(K)\subseteq V\setminus V_i$ be the set $\{g\cdot w_i: g\in H_K\}$.  In a similar way to one above, we have that
\[
|V(K)|=\sum_{w\in V(K)}\frac{|\{g\in H_K:g\cdot w_i=w\}|}{|\mathrm{Stab}_{w}\cap H_K|}.
\]
Therefore, we conclude that for any $\emptyset \ne K \subseteq V$,
\begin{align*}
|V(K)|&\geq \sum_{w\in V(K)}\frac{|\{g\in H_K:g\cdot w_i=w\}|}{|\mathrm{Stab}_{w}|}\\
&=\frac{\sum_{w\in V(K)}|\{g\in H_K:g\cdot w_i=w\}|}{|\mathrm{Stab}_{v_i}|}=\frac{|H_K|}{|\mathrm{Stab}_{v_i}|}=|K|.
\end{align*}
The marriage theorem, therefore, verifies our claim (note that $V(K)$ coincides with the set $\bigcup_{v\in K}\{g\cdot w_i:g\in G_{i,i}, \ g\cdot v_i=v\}$).

Fix $i\in [n]$. Set  $A_i^{(1)}=A_i\cap V_i$ and $A_i^{(2)}=A_i-V_i$. Note that by item $(c)$ in Theorem~\ref{mthm:imprimitive}, $|A_i^{(2)}|\leq \frac{4h_n(\Gamma)}{h_{n+1}(\Gamma)}|V|$. Then the claim above implies that 
\[
|\partial (A_i^{(1)}, V\setminus A_i) |\geq |\partial (A_i^{(1)}, V\setminus (V_i\cup A_i^{(2)}) )| \geq |A_i|- \frac{8h_n(\Gamma)}{h_{n+1}(\Gamma)}|V|.
\] 
We hence have that 
\[
\frac{|\partial A_i|}{|A_i|}\geq 1 -\frac{8h_n(\Gamma)}{h_{n+1}(\Gamma)}\frac{|V|}{|A_i|}.
\]
Take the minimum over all $i\in [n]$. Then, by definition, the minimum of the left-hand side equals $h_n(\Gamma)$. By $(\diamond)$ in Section~\ref{section:Phi}, we conclude that
\[
h_n(\Gamma)\geq 1- \frac{n^2+1}{n}\frac{8h_n(\Gamma)}{h_{n+1}(\Gamma)}\geq 1- 10n\cdot \frac{h_n(\Gamma)}{h_{n+1}(\Gamma)}
\]
because $n\geq 2$. These inequalities lead us to  the desired inequality on $h_n$'s. 

For the inequalities on $\iota_n$'s, in a similar manner to the one above, we can show that for every $i\in [n]$, 
\[
\frac{|\delta A_i|}{|A_i|}\geq 2 -\frac{16\iota_n(\Gamma)}{\iota_{n+1}(\Gamma)}\frac{|V|}{|A_i|}.
\]
This ends our proof of Theorem~\ref{mthm:mt}.
\end{proof}
\section{Remark on vertex and edge transitive graphs}\label{section:ve}

Theorem~\ref{mthm:mt} shows that, if we restrict the class of regular finite connected graphs to that of \textit{vertex-transitive} graphs, then we obtain a non-trivial inequality between $h_{n+1}(\Gamma)$ and $h_n(\Gamma)$. In this class, as we saw in Section~\ref{section:counterexample}, item $(1)$ of Question~\ref{question:M} resolves in the negative. Here, we mention that, if we stick to a much smaller class, that of \textit{vertex and edge transitive} graphs, then this question has the positive answer.

\begin{corollary}\label{cor:vertexandedge}
Let $\Gamma$ be a finite connected graph. If $\Gamma$ is vertex and edge transitive, then for any $2\leq n \leq |V|-1$, we have that $h_{n+1}(\Gamma)\leq (10n+1)h_{n}(\Gamma)$.
\end{corollary}

\begin{proof}
Suppose that $h_{n+1}(\Gamma)> 2(n+1)h_{n}$. Then by Theorem~\ref{mthm:imprimitive}, there exists a system $(V_1,\ldots ,V_n)$ of imprimitivity of size $n$ for $\mathrm{Aut}(\Gamma)\curvearrowright \Gamma$. If there exists an edge inside $V_i$ for some $i$, then  it contradicts the assumption. Indeed, since $\Gamma$ is connected and the group action is vertex-transitive, then there must exist $v,v'\in V_i$ and $w\in V\setminus V_i$ such that $(v,v')$ and $(v,w)$ are in $E$. By the edge-transitivity, this contradicts the imprimitivity of the system. 

There are, hence, no edges inside $V_i$ for each $i$. Then by item $(c)$ of Theorem~\ref{mthm:imprimitive}, in a similar argument to one in Section~\ref{section:mt}, we have that 
\[
h_n(\Gamma)\geq d_{\Gamma}- d_{\Gamma}\cdot \frac{n^2+1}{n}\frac{8h_n(\Gamma)}{h_{n+1}(\Gamma)}.
\]
This implies that 
\[
h_n(\Gamma)\geq \frac{d_{\Gamma}h_{n+1}(\Gamma)}{10d_{\Gamma}n+h_{n+1}(\Gamma)}.
\]
Because $h_{n+1}(\Gamma)\leq d_{\Gamma}$, we obtain the conclusion.
\end{proof}

\section{Further questions}\label{section:que}
\begin{question}
\begin{enumerate}[$(1)$]
  \item For $n\geq 2$, construct an example of a finite connected vertex-transitive graph $\Gamma$ such that $h_{n+1}(\Gamma) >2(n+1)h_n(\Gamma)$ and any realizer of $h_n(\Gamma)$ is not a system of imprimitivity of size $n$ for $\mathrm{Aut}(\Gamma) \curvearrowright V$.
  \item If we know that a vertex-transitive finite graph $\Gamma$ does not admit a system of imprimitivity of size $n$ for $\mathrm{Aut}(\Gamma) \curvearrowright V$, then $h_{n+1}(\Gamma)\leq 2(n+1)h_n(\Gamma)$ $($see Section~$\ref{section:npoints}$$)$. Is it possible to find a universal constant $C>0$, even idependent of $n$, such that for any $n<|V|$, $h_{n+1}(\Gamma)\leq Ch_n(\Gamma)$ holds in the setting above?
  \item Establish a similar inequality, to ones in Theorem~$\ref{mthm:mt}$, between $\lambda_{n+1}(\Gamma)$ and $\lambda_n(\Gamma)$. Note that, as we saw in Section~$\ref{section:counterexample}$, there must exist a critical value of $\lambda_n(\Gamma)$ to bound $\lambda_{n+1}(\Gamma)$ from above.
   \item For a finite connected vertex-transitive  graph, obtain a direct comparison between $\iota_n(\Gamma)$ and $\iota_2(\Gamma)$. Similar questions on $h_n$ and on $\lambda_n$, respectively, may be considered. However, in these cases, we need to impose appropriate conditions on $\Gamma$ to avoid examples in Section~$\ref{section:counterexample}$.
\end{enumerate}
\end{question}

\section*{acknowledgments}
The author is grateful to Koji Fujiwara for asking his question to the author, and for explaining his idea. The author thanks Kei Funano, Masaki Izumi, Akihiro Munemasa, Hidehiro Shinohara, and Mamoru Tanaka for helpful comments and references. He gratefully acknowledges the two referees for their perserving reading of the submitted version of the current manuscript; for various comments and suggestions, that considerably improve the present paper; and for providing the author with relevant references.

\bibliographystyle{amsalpha}
\bibliography{mimura.bib}

\end{document}